\documentclass[reqno,12pt]{amsart}      
\usepackage{amsfonts}
\usepackage{colordvi}
\usepackage{amsmath,amssymb,amsthm} 
\usepackage{mathrsfs}
\usepackage{graphicx}
\usepackage{ulem}  
\theoremstyle{plain}

\theoremstyle{remark}

\newcounter{FNC}[page]
\def\fauxfootnote#1{{\addtocounter{FNC}{2}\Magenta{$^\fnsymbol{FNC}$}%
     \let\thefootnote\relax\footnotetext{\Magenta{$^\fnsymbol{FNC}$#1}}}}

\newcommand{\defcolor}[1]{\Blue{#1}}
\newcommand{\demph}[1]{\defcolor{{\sl #1}}}

\newcommand{\Gr}{\mbox{\rm Gr}}
\newcommand{\CC}{{\mathbb C}}
\newcommand{\PP}{{\mathbb P}}

\newcommand{\alphadot}{\alpha^\bullet}

\begin{document}     


\title[Numerical Schubert Calculus in Macaulay2]{Numerical Schubert Calculus in Macaulay2}

\author[Leykin]{Anton Leykin}
\address{School of Mathematics\\
         Georgia Institute of Technology\\
         686 Cherry Street\\
         Atlanta, GA 30332-0160\\
         USA}
\email{leykin@math.gatech.edu}
\urladdr{http://people.math.gatech.edu/~aleykin3/}
\author[Mart\'in del Campo]{Abraham Mart\'in del Campo}
\address{Abraham Mart\'in del Campo\\
         Centro de Investigaci\'on en Matem\'aticas, A.C.\\
         Jalisco S/N, Col. Valenciana\\
         Guanajuato, Gto.  M\'exico}
\email{abraham.mc@cimat.mx}
\urladdr{http://personal.cimat.mx:8181/~abraham.mc/Home.html}
\author[Sottile]{Frank Sottile}
\address{Frank Sottile\\
         Department of Mathematics\\
         Texas A\&M University\\
         College Station\\
         Texas \ 77843\\
         USA}
\email{sottile@math.tamu.edu}
\urladdr{www.math.tamu.edu/\~{}sottile}
\author[Vakil]{Ravi Vakil}
\address{Ravi Vakil\\
      Department of Mathematics\\ 
      Stanford University\\ 
      Stanford, CA 94305
       USA}
\email{rvakil@stanford.edu} 
\urladdr{http://math.stanford.edu/\~{}vakil}
\author[Verschelde]{Jan Verschelde}
\address{Jan Verschelde\\
      Dept of Math, Stat, and CS\\
      University of Illinois at Chicago\\ 
      851 South Morgan (M/C 249)\\ 
      Chicago, IL 60607
      USA} 
\email{jan@math.uic.edu} 
\urladdr{http://www.math.uic.edu/\~{}jan}
\thanks{The authors thank the American Institute of Mathematics for supporting the project
  through their SQuaREs program.}
\thanks{The package {\tt NumericalSchubertCalculus} is included in Macaulay2 version 1.18.}

\begin{abstract}
 The Macaulay2 package {\tt NumericalSchubertCalculus} provides methods for the numerical computation of
 Schubert problems on Grassmannians.
 It implements both the Pieri homotopy algorithm and the Littlewood-Richardson homotopy algorithm.
 Each algorithm has two independent implementations in this package.
 One is in the scripting language of Macaulay2 using the 
 package {\tt NumericalAlgebraicGeometry}, and the other is in the compiled code 
 of PHCpack.
\end{abstract}

\maketitle

\section{Introduction}
The Schubert calculus on the Grassmannian involves all problems of determining the linear subspaces of a vector space that have 
specified positions with respect to fixed flags of linear subspaces.
The enumeration of the solution linear spaces may be solved using the Macaulay2 package {\tt Schubert2}.
Numerical Schubert calculus computes the actual solution planes to a given instance of a problem from the 
Schubert calculus using numerical methods, and the eponymous Macaulay2 package implements algorithms that
accomplish this task. 

Schubert problems arise in applications in control theory~\cite{By89} and in information theory~\cite{BCT}, and
they form a rich class of geometric problems that serves as a laboratory for investigating new phenomena in
enumerative geometry such as reality~\cite{RSEG} and Galois groups~\cite{MS}.
The ability to compute solutions has been important in these areas.
Schubert problems are also challenging to solve using standard numerical methods.
This is because Schubert problems are typically not complete intersections, and even when they are, they have far fewer solutions than
standard combinatorial bounds~\cite[p.~768]{HSS98}.

{\tt NumericalSchubertCalculus} has methods implementing numerical homotopy continuation 
algorithms that exploit explicit geometric proofs of multiplication formulas in the cohomology of a Grassmannian.
The Pieri homotopy algorithm~\cite{HSS98} is based on the geometric Pieri formula~\cite{So97} and the
Littlewood-Richardson homotopy algorithm~\cite{LRH,SVV} is based on the geometric Littlewood-Richardson
rule~\cite{Va06a}. 
{\tt NumericalSchubertCalculus} has two implementations of each algorithm.
One is in the Macaulay2 scripting language using its {\tt NumericalAlgebraicGeometry}~\cite{NAG4M2} package
and the other is compiled code using PHCpack~\cite{PHCpack}.

\section{Mathematical background}\label{S:background}

A Schubert problem on a Grassmannian is a problem of determining the linear subspaces of
a given dimension (its solutions) that have specified positions with respect to other
fixed, but general, linear subspaces, when there are finitely many solutions. 
The simplest non-trivial Schubert problem asks for the two-dimensional
subspaces $H$ of $\CC^4$ that have nontrivial intersection with each of four general two-dimensional linear
subspaces $L_1,\dotsc,L_4$. 
Replacing $\CC^4$ by projective space, this becomes the problem of determining the lines $h$ in $\PP^3$ that
meet four general lines $\ell_1,\dotsc,\ell_4$.

\begin{figure}[htb]
%
%
\centerline{
  \begin{picture}(189,108)(-1,5)
   \put(3,0){\includegraphics[height=4.1cm]{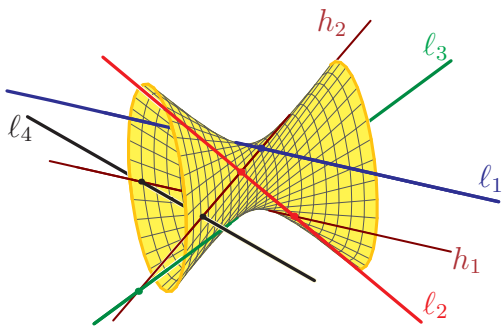}}
   \put(159,  5){\Red{$\ell_2$}}
   \put(180, 53){\Blue{$\ell_1$}}
   \put(159,103){\Green{$\ell_3$}}
   \put(  2, 72){$\ell_4$}
   \put(170, 22){\Maroon{$h_1$}}
   \put(119,112){\Maroon{$h_2$}}
  \end{picture}
}
 \caption{Problem of four lines}\label{F:4lines}
\end{figure}

This problem has two solutions.
To see this, we use the classical observation that three mutually skew lines 
\Blue{$\ell_1$}, \Red{$\ell_2$}, and \Green{$\ell_3$} 
lie on a unique hyperboloid (Fig.~\ref{F:4lines}).
This hyperboloid has two rulings, one contains 
\Blue{$\ell_1$}, \Red{$\ell_2$}, and \Green{$\ell_3$}, 
and the second consists of the lines meeting these three.
If the fourth line, $\ell_4$, is general, then it will meet the hyperboloid in two points,
and through each of these points there is a unique line in the second ruling.
These two lines, \Maroon{$h_1$} and \Maroon{$h_2$}, are the solutions to this instance
of the problem of four lines.

Let \defcolor{$\Gr(k,n)$} be the Grassmannian of $k$-dimensional linear subspaces ($k$-planes) in $\CC^n$.
This has dimension $k(n{-}k)$.
Indeed, a general linear subspace $H\in\Gr(k,n)$ is the column space of a matrix 
$\left(\begin{smallmatrix}X\\I_k\end{smallmatrix}\right)$ in (reverse) column-reduced echelon form, where $X$ is a
$(n{-}k)\times k$ matrix, and different matrices determines a different $k$-planes.
An incidence condition on $k$-planes is encoded by a \demph{bracket}, which is an
increasing sequence $\defcolor{\alpha}\colon 1\leq\alpha_1<\alpha_2<\dotsb<\alpha_k\leq n$ of integers.
The condition is imposed by a flag \defcolor{$F$} of linear spaces 
$F\colon F_1\subset F_2\subset\dotsb\subset F_n=\CC^n$ where $\dim F_i=i$.
This pair defines a \demph{Schubert variety},
 \begin{equation}\label{Eq:SchubertVariety}
   \defcolor{X_\alpha F}\ :=\ \{ H\in\Gr(k,n)\mid \dim H\cap F_{\alpha_i}\geq i\,\ \mbox{for}\  i=1,\dotsc,k\}\,.
 \end{equation}
Requiring that $H\in X_\alpha F$ is a \demph{Schubert condition} on $H$ of type $\alpha$ imposed by the flag $F$.

The Schubert variety $X_\alpha F$ has dimension $|\alpha|:=\sum_i \alpha_i-i$ and thus codimension 
$\defcolor{\|\alpha\|}:=k(n{-}k)-|\alpha|$.
A \demph{Schubert problem $\alphadot$} is a list $\alpha^1,\dotsc,\alpha^s$ of brackets that satisfy
$\sum_i\|\alpha^i\|=k(n{-}k)$.
An \demph{instance} of $\alphadot$ is given by a list $F^1,\dotsc,F^s$ of flags, and is the geometric problem of
those $H\in\Gr(k,n)$ that satisfy the Schubert condition $\alpha^i$ imposed by $F^i$ ($H\in X_{\alpha^i}F^i$) for each $i$.
These form the intersection
 \begin{equation}\label{Eq:intersection}
  X_{\alpha^1}F^1\,\bigcap\,
  X_{\alpha^2}F^2\,\bigcap\, \dotsb \,\bigcap\,
  X_{\alpha^s}F^s\,.
 \end{equation}
For an example, consider the problem of four lines (expressed in $\Gr(2,4)$), let $F^i$ be a flag with $2$-plane 
$F^i_2=L_i$, for each $i=1,\dotsc,4$.
The corresponding Schubert condition is that $\dim H\cap F^i_2\geq 1$ and $\dim H\cap F^i_4=2$ (as
$F^i_4=\CC^4$), and is given by the bracket $\{2,4\}$, and so  the problem of four lines is
$\alphadot=(24,24,24,24)$ (write $24$ for the bracket $\{2,4\}$).

Kleiman~\cite{Kl74} proved that if the flags are general, then the intersection~\eqref{Eq:intersection} is transverse and 
there are finitely many solutions.
The package {\tt NumericalSchubertCalculus} implements methods to compute the solutions~\eqref{Eq:intersection} to a given instance of a
Schubert problem.

The number \defcolor{$d(\alphadot)$} of solutions (the points in~\eqref{Eq:intersection}) is independent of choice of general flags,
and this may be computed using algorithms in the Schubert calculus.
These are implemented in the {\tt Macaulay2} package {\tt Schubert2}.
A geometric derivation of $d(\alphadot)$ was given by the geometric Littlewood-Richardson rule in~\cite{Va06a}.
The geometric deformations of the geometric Littlewood-Richardson rule underlie the Littlewood-Richardson homotopy
algorithm~\cite{LRH,SVV} which is implemented in {\tt NumericalSchubertCalculus}.

A Schubert condition $\alpha$ is \demph{simple} if $\|\alpha\|=1$.
A Schubert problem in which all conditions except possibly two are simple is a \demph{simple Schubert problem}.
A geometric proof of the Pieri rule~\cite{So97} led to the Pieri homotopy algorithm~\cite{HSS98} for solving simple Schubert problems, and
this is also implemented in {\tt NumericalSchubertCalculus}.

\section{Implementation and Syntax}

The package {\tt NumericalSchubertCalculus} has methods to compute the solutions to a given instance~\eqref{Eq:intersection} of a Schubert
problem.
It includes two independent implementations of both the Littlewood-Richardson homotopy
algorithm~\cite{LRH,SVV}  and the Pieri homotopy algorithm~\cite{HSS98}.
One is in the Macaulay2 scripting language and the other is compiled code.
These implementations have slightly different capabilities and input syntax, with the package providing some interoperability.
We briefly describe the input and output syntax of the primary methods, and their capabilities.

In both, brackets are represented by lists so that {\tt \{3, 5, 6\}} is a simple
Schubert condition when $k=3$ and $n=6$.
Flags are represented by invertible $n\times n$ matrices, with $F_i$ the span of the first $i$ columns.
An element $H$ of $\Gr(k,n)$ is represented by an $n\times k$ matrix whose column span is $H$.

\subsection{Scripted methods}\label{SS:NSC_M2}
These compute the solutions to a given instance of a Schubert problem.
An instance of a Schubert problem $\mbox{\tt SchubProb}=\{\alpha^1,\dotsc,\alpha^s\}$ on $\Gr(k,n)$ is represented by a list of pairs
 \begin{equation}\label{Eq:SchProbInst}
   \mbox{\tt SchProbInst}\ =\
     \{\, \{\alpha^1,F^1\}\,,\, \{\alpha^2,F^2\}\,,\dotsc\,,\, \{\alpha^s,F^s\}\,\}\,,
 \end{equation}
where each pair $\{\alpha^i,F^i\}$ consists of a bracket $\alpha^i$ and a flag $F^i$.
The parameters $k,n$ are implicit in the data ($k$ is the length of the bracket and $n$ is the size of the matrix).

Given an instance~\eqref{Eq:SchProbInst} of a Schubert problem on $\Gr(k,n)$, the solutions  are computed with the method
{\tt solveSchubertProblem}.
A typical call is
\[
   \mbox{\tt Solns\ =\ solveSchubertProblem(SchProbInst, k, n)}.
\]
(The parameters $k,n$ are included for internal verification that {\tt SchProbInst} is in fact an instance of a Schubert problem on
$\Gr(k,n)$.)
This returns the solutions as a list {\tt Solns} of $n\times k$ matrices $H_1,\dotsc,H_d$, each of whose column span is a $k$-plane solving
the instance represented by {\tt SchProbInst} that is, it is a point in the intersection~\eqref{Eq:intersection}.
As the computation operates in local coordinates for the Grassmannian, it presupposes that the flags are sufficiently general so that the
solutions lie in the local coordinates.

The method {\tt randomSchubertProblemInstance} returns a random instance of a Schubert problem
$\mbox{\tt SchubProb}=\{\alpha^1,\dotsc,\alpha^s\}$ on $\Gr(k,n)$,  
\[
   \mbox{\tt RandSchProbInst\ =\ randomSchubertProblemInstance(SchubProb, k, n)}\,.
\]
This gives an instance~\eqref{Eq:SchProbInst} of {\tt SchubProb} in which the flags $F^i$ are random complex matrices.
Using the output {\tt RandSchProbInst} as input for {\tt solveSchubertProblem} gives solutions to the random instance, and not a user-chosen
instance.
{\tt NumericalSchubertCalculus} provides a method, {\tt changeFlags}, based on the parameter or cheater's homotopy~\cite{LSY89,MS89} that,
given the solutions {\tt Solns} to a particular instance of a Schubert problem {\tt SchubProb}, represented by a list of flags {\tt F},
computes the solutions to a different instance represented by another list of flags {\tt myF}. 
A typical call is
\[
\mbox{\tt mySols\ =\ changeFlags(Solns, CFG)}\,,
\]
where {\tt CFG} is a triple  {\tt \{SchubProb, F,  myF\}}

If {\tt SchubProb} is a simple Schubert problem, then the alternative Pieri homotopy algorithm is available, and it may be used to
compute the solutions with the call 
\[
   \mbox{\tt Solns\ =\ solveSimpleSchubert(SchubProb, k, n)}\,.
\]
This requires that all conditions in the instance {\tt SchubProb} are simple, except possibly the first two.
The previous discussion of a random versus a user-provided instance and the method {\tt changeFlags} for interpolating between them also
applies  here.

One exported method is {\tt checkIncidenceSolution}, which is a Boolean-valued function that may be used to check if a given $k$-plane $H$
satisfies an instance of a Schubert problem.
A typical call is
\[
   \mbox{\tt checkIncidenceSolution(H, SchProbInst)}\,.
\]
Here, {\tt H} is a $n\times k$ matrix and {\tt SchProbInst} is a list of pairs~\eqref{Eq:SchProbInst}.

The exported method {\tt LRnumber} computes the number of solutions to a given Schubert problem {\tt SchubProb}.
A typical call is
\[
   \mbox{\tt d\ =\ LRnumber(SchubProb, k, n)}\,,
\]
where {\tt SchubProb} is a list of brackets that constitute a Schubert problem on $\Gr(k,n)$ and the output {\tt d} is an integer.

While it is straightforward~\eqref{Eq:SchubertVariety} to encode a Schubert condition in $\Gr(k,n)$ as a bracket
$\alpha\colon 1\leq\alpha_1<\dotsb<\alpha_k\leq n$,
a common equivalent encoding is a partition, which is a weakly decreasing sequence
$\lambda\colon n{-}k\geq\lambda_1\geq\dotsb\geq\lambda_k\geq 0$.
For example, {\tt Schubert2} uses partitions.
A bracket $\alpha$ and its corresponding partition $\lambda$ satisfy
\[
   \alpha_i{-}i\ +\ \lambda_i\ =\ n{-}k\quad\mbox{for }\  i=1,\dotsc,k\,.
\]
In this case, $|\lambda|=\lambda_1+\dotsb+\lambda_k$ equals the codimension $\|\alpha\|$.
The exported methods {\tt bracket2partition} and {\tt partition2bracket} translate between the two notations.
The input for methods described so far allow a Schubert problem to be expressed as either a list of brackets or a list
of partitions.

\subsection{Interface to PHCpack}\label{SS:PHC}
The corresponding routines that interface with the compiled implementations in PHCpack have minor differences in syntax and capability to
those in Section~\ref{SS:NSC_M2}.
The main difference is that these methods compute solutions to a random instance of the given Schubert problem that is generated when the
method is called.
The output includes the instance that was used for the computation.
To obtain solutions to a user's instance requires using {\tt changeFlags}. 

The method {\tt LRtriple} implements the geometric Littlewood-Richardson rule.
It encodes a Schubert problem slightly differently than the method {\tt solveSchubertProblem};
rather than a list of brackets, it takes a matrix whose rows have the form {\tt \{m, b\}} where
{\tt b} is a bracket (increasing sequence of length $k$) and {\tt m} is the \demph{multiplicity}
of that bracket in the Schubert problem (how often it appears).
The method {\tt NSC2phc} takes as input a list of brackets and its output is a matrix encoding the Schubert problem.
A typical call is
\[
  {\tt M\ =\ NSC2phc(SchubProb, k, n))}\,.
\]
Given a matrix {\tt M} encoding a Schubert problem on $\Gr(k,n)$, a typical call of   {\tt LRtriple} is
\[
   {\tt (F, P, S)\ =\    LRtriple(n, M)}\,.
\]
The output is a triple of strings with {\tt F} describing the flags and the solutions, {\tt P} is the polynomial system solved in local
coordinates, and {\tt S} are the solutions in local coordinates and information about the computation (see the documentation).
The method {\tt parseTriplet} transforms these strings into Macaulay2 objects for further processing.
A typical call is
\[
  {\tt (R, pols, sols, fixedFlags, movedFlag, solutionPlanes)\ =\ parseTriplet(F,P,S) }\,.
\]
This is a list of Macaulay2 objects which enable further processing of the computed solutions (again, see the documentation).


The routine {\tt LRrule} is a PHCpack implementation of a method to compute the number of solutions to a given Schubert problem using
the geometric Littlewood-Richardson rule.
A typical call is
\[
    {\tt s\ =\  LRrule(n,M)}\,,
\]
where {\tt n} is the ambient dimension, and {\tt M}  is a matrix encoding the Schubert problem.
The output {\tt s} is a string encoding a product in the cohomology ring.
For problem of four lines, this is {\tt [ 2 4 ]\^{}4 = +2[1 2]}.

\subsection{A note on numerics}
{\tt NumericalSchubertCalculus} implements highly recursive algorithms based on numerical homotopy continuation and
path-tracking.
Because of the recursion, when running the scripted algorithms of Section~\ref{SS:NSC_M2}, it is advisable to increase
{\tt Macaulay2}'s limit on the allowed depth of a recursion, which is controlled by the {\tt recursionLimit} variable.
As the algorithms are numerical, at times the path-tracking will fail.
One work-around is to set the seed for the random number generator.
For the scripted algorithms of Section~\ref{SS:NSC_M2}, this may be done by changing the variable {\tt setRandomSeed}, and for the
compiled algorithms in Section~\ref{SS:PHC}, this is done with the option {\tt   RandomSeed}.
The numerical precision may be set in the compiled algorithms in Section~\ref{SS:PHC} using the option {\tt WorkingPrecision}.

Numerical computation in {\tt NumericalSchubertCalculus} takes place in systems of local Stiefel coordinates for subvarieties of $\Gr(k,n)$.
While details are explained in~\cite[\S2.2]{LRH}, we note the following consequence:
When computing solutions to a Schubert problem $\alpha^1,\dotsc,\alpha^s$, the maximum number of coordinates (the ambient dimension of the
computation) is $k(n{-}k)-\|\alpha^1\|-\|\alpha^2\|$.
The computational heuristic that using fewer coordinates improves performance holds in practice with our software.
Both our implementations of the Littlewood-Richardson homotopy, {\tt solveSchubertProblem} and {\tt LRtriple}, are faster and more
numerically stable when the conditions satisfy
\[
\min\{ \|\alpha^1\|,\|\alpha^2\|\}\ \geq\ \max\{\|\alpha^3\|,\dotsc,\|\alpha^s\}\}\,.
\]
Consequently, a user is advised to sort their Schubert conditions so that the largest two are first.

\section{Examples}

We close with some examples that illustrate our software.
The problem of four lines in $\Gr(2,4)$ is given by four brackets $\{\{2,4\},\{2,4\},\{2,4\},\{2,4\}\}$.
{\small
\begin{verbatim}
i4 : SchubProb = randomSchubertProblemInstance({{2,4},{2,4},{2,4},{2,4}}, 2,4)

o4 = {({2, 4}, | -.114281-.993448ii -.0168841-.999857ii -.0141397-.9999ii....
               | -.788054+.615606ii -.319468+.947597ii  .980731-.195363ii....
               | .300518-.953776ii  .0796792-.996821ii  .624656-.7809ii  ....
               | .676179+.736737ii  -.886912+.461938ii  -.237463+.971397i....
     --------------------------------------------------------------------....

o4 : List

i5 : solveSchubertProblem(SchubProb,2,4)

o5 = {| -.140237+1.14381ii -.378754+.204742ii |, | .0466437-.65223ii    -....
      | 1.57394-.065106ii  .0818881+.152422ii |  | -.0559878+.7121ii    1....
      | -.906612+.910475ii .0437461-.383392ii |  | -.00398277-.653482ii -....
      | -2.04662-1.34516ii .0997269+.124121ii |  | -1.14257+.21797ii    1....

o5 : List
\end{verbatim}
}
\noindent
The output {\tt o5} consists of two $4\times 2$ matrices whose column spaces solve the given instance of this Schubert
  problem.
Its output is nearly instantaneous.
Here is a significantly larger Schubert problem in $\Gr(4,8)$ with 1530 solutions, computed using {\tt LRtriple}:
{\small
\begin{verbatim}
i1 : M = matrix {{2, 3,5,7,8 },{1, 3,6,7,8}, {8, 4,6,7,8}};

               3        5
o1 : Matrix ZZ  <--- ZZ

i2 : result = LRtriple(8, M);

i3 : L= lines(result_2);

i4 : L_1

o4 = 1530 10

i5 : L_(#L-10)

o5 = User time in seconds was                    336.834157000 =  0h 5m36s834ms
\end{verbatim}
}
We now compute the simple Schubert problem $\{3,5,6\}^9$ in $\Gr(3,6)$ with 42 solutions in four different ways.
  The bracket $\{3,5,6\}$ corresponds to the partition $\{1\}$.
{\small
\begin{verbatim}
i4 : k=3, n=6;

i5 : conds = {{1},{1},{1},{1},{1},{1},{1},{1},{1}};

i6 : LRnumber(conds,k,n)

o6 = 42

i7 : SchPblm = randomSchubertProblemInstance(conds,k,n);

i8 : time S = solveSchubertProblem(SchPblm,k,n);
     -- used 96.698 seconds

i9 : #S

o9 = 42

i10 : time S = solveSimpleSchubert(SchPblm,k,n);
     -- used 0.877903 seconds

i11 : #S

o11 = 42

i12 : M = NSC2phc(conds,k,n)

o12 = | 9 3 5 6 |

               1        4
o12 : Matrix ZZ  <--- ZZ

i13 : s = LRrule(n,M)

o13 = [ 3 5 6 ]^9 = +42[1 2 3]

i14 : time (F,P,S) = LRtriple(n,M);
     -- used 0.122594 seconds

i15 : (lines(S))_1

o15 = 42 7

i16 : time (ipt, otp) = PieriHomotopies(3,3);
     -- used 0.0915026 seconds

i17 : ipt

o17 = {| -.408248           .233797-.146195ii   -.339562+.0273021ii  |, ....
       | .403083-.0647356ii -.361496+.0525928ii .170936-.210227ii    |  ....
       | -.35011+.209975ii  .181004-.236536ii   -.245246+.0877148ii  |  ....
       | -.332084+.237459ii -.232994+.556345ii  .121023+.101895ii    |  ....
       | .403492+.0621347ii .524032-.051926ii   .0282015-.00137244ii |  ....
       | .168565-.371823ii  .119146+.21564ii    -.0367667+.845884ii  |  ....

o17 : List

i18 : #otp

o18 = 42
\end{verbatim}}


\providecommand{\bysame}{\leavevmode\hbox to3em{\hrulefill}\thinspace}
\providecommand{\MR}{\relax\ifhmode\unskip\space\fi MR }
\providecommand{\MRhref}[2]{%
  \href{http://www.ams.org/mathscinet-getitem?mr=#1}{#2}
}
\providecommand{\href}[2]{#2}

\end{document}